\documentclass{article}
\usepackage{amsfonts,amsmath}
\usepackage{amssymb}
\usepackage{amsthm}
\usepackage[cp1251]{inputenc}

\newtheorem{thm}{Theorem}[section]

\newtheorem{cor}[thm]{Corollary}
\newtheorem{proposition}[thm]{Proposition}

\newtheorem{rem}[thm]{Remark}

\begin{document}

\title
{Every bounded self-ajoint operator is a real linear combination
of $4$ orthoprojections}

\author{V. Rabanovich}

\maketitle
\begin{abstract}
We prove that every bounded self-adjoint operator in Hilbert space
is a real linear combination of $4$ orthoprojections. Also we show
that operators of the form identity minus compact positive
operator can not be decomposed in a real linear combination of $3$
orthoprojections. Using ideas applied in infinite dimensional
space, we find $n\times n$ matrices that are not real linear
combinations of $3$ orthoprojections for every $n\ge 76$.
\end{abstract}



\sloppy

\newcommand{\tr}{\mathrm{tr}\,}
\newcommand{\rank}{\mathrm{rank}\,}
\newcommand{\diag}[1]{\mathrm{diag}\,(#1)}
\newcommand{\Diag}[1]{\mathrm{diag}\,\left(#1\right)}
\renewcommand{\Im}{\mathrm{Im}\,}
\newcommand{\Ker}{\mathrm{Ker}\,}




\section{Introduction}
We consider here linear combinations of orthogonal projections
$P_i$, $P_i^*=P_i^2=P_i$, on a separable infinite-dimensional
Hilbert space $H$. It was proved in \cite{Fillmore67} that every
bounded  operator on $H$ is a complex linear combination of $257$
orhthoprojections. Later the author showed that every bounded
self-adjoint operator $A$ on $H$ is a real linear combination of
$9$ orhthoprojections \cite{Fill69}. At the same time C.~Pearcy
and D.~Topping proved that it is enough only $8$ summands in the
real combination \cite{PeTo67}. Then
 it was established in \cite{Paszkiew80} that
 the operator $A$  is a real linear combination of $6$ orhthoprojections.
By modification of the proof in \cite{Paszkiew80}, K.~Matsumoto
diminished the number of summands in such a decomposition to $5$
items \cite{Matsumoto84}. The fact that not every Hermitian
operator can be decomposed into a linear combination of two
orthoprojections is known for many years due to some symmetry
property of a linear combination of two orthoprojections (see
Proposition \ref{prop1}). The most simple example of such an
operator is a Hermitian operator in $4$ dimensional space with
eigenvalues $0.9$, $1$, $1.01$ and $1.0001$.  So it remains to
consider real combinations of $3$ or $4$ orthprojections. It was
proved in \cite{RVI0305} that every diagonizable self-adjoint
operators is a real linear combination of $4$ orthoprojections.
Using some simple manipulations with self-commutators in Theorem
\ref{th1}, we give a short possibly new proof of this result for
every self-adjoint operator. It should be noted that the authors
from \cite{GolPasz92} considered decompositions of operators into
a linear combinations of $4$ orthoprojections in von Neumann
algebras, factors of type $I$. They formulated the same result as
Theorem \ref{th1} but with reference in  their proof to
unpublished paper.

 As a simple corollary of Theorem
\ref{th1}, we find that self-adjoint operator can be decomposed
into an integral combination of $5$ orthoprojections (the fact
 that was prove in \cite{GolPasz92}[Theorem 1(3b)] also), and into a
real combination of $5$ orthoprojections with infinite dimensional
kernel and range.

In section \ref{sec4} we give a class of operators for which there
is no decompositions into a linear combinations of $3$
orthoprojections. It appears that these are the operators of the
form $I+K$, where $I$ is the identity operator and $K$ is an
infinite-rank compact negative or positive operator.

It directly follows from our result that every bounded operator is
a complex linear combination of $8$ orthoprojections. We can not
show that such a number is minimal. Instead of this we give in
Corollary \ref{cor5} an example of operator which is not a complex
linear combination $4$ orthoprojections.

At the end of the paper we consider finite-dimensional unitary
space. For a Hermitian $n\times n$ matrix $A$, Y.~Nakamura
\cite{Nakamura84} proved that $A$ is  a linear combination of $4$
orthoprojections, and it is a linear combination of $3$
orthoprojections for $n\le 7$. Using ideas from infinite
dimensional case, we find $m\times m$ matrices $B_m$, $m \ge 76$,
a small norm perturbation of a scalar matrix, such that $B_m$ can
not be presented as a linear combinations of $3$ orthoprojection.


Throughout the paper $X\approx Y$ means that $X$ is similar to $Y$
and $\diag{a_1,\dots,a_n}$ means a \emph{diagonal} or \emph{block
diagonal} matrix with diagonal elements $a_1,\dots,a_n$ from
$\mathbb{C}$ or from algebra of bounded operator on a Hilbert
space $H$. We denote by $\tr A$ the \emph{trace} of $A$ and by
$\sigma(A)$ its \emph{spectrum}. All eigenvalues
$\lambda_1(A),\dots, \lambda_n(A)$ counting multiplicity of a
Hermitian $n\times n $ matrix $A$ will supposed to be arranged in
increasing order, $\lambda _1(A)\le \lambda_2(A)\le \dots\le
\lambda_n(A)$. Also we set $\mathbf{0}_n=\diag{0,0,\dots,0}$, and
$I_n=\diag{1,1,\dots,1}$. Identity operator on $H$ will be denoted
by $I$ or $I_H$ and zero operator --- by $0_H$.

\section{Preliminaries}
Before we start with linear combinations of $4$ orthoprojections,
we remind some facts on linear combinations of two
orthoprojections. Suppose $P_1$ and $P_2$ are orthoprojections on
a Hilbert space $H$. If $v\in H$ is  an eigenvector of both $P_1$
and $P_2$, then for every $a,b\in \mathbb{R}$, it is an
eigenvector of the operator $A=aP_1+bP_2$ with eigenvalue $\lambda
\in \{0,a,b,a+b\}$. It follows from \cite{Halmos2s}, that the
inverse statement is also true, that is for $\mu \in
\{0,a,b,a+b\}$ and $h\in H$, the equality $Ah=\mu h$ yields $h$ is
an eigenvector of both operators $P_1$ and $P_2$. Also every point
$x$ from $\sigma(A)$ lies in the union of two segments $[0,a]$ and
$[b,a+b]$ for $|b|\ge |a|$:
\begin{equation}\label{eq2}
x\in [0,a]\cup[b,a+b]
\end{equation}
  The following Proposition is a direct
corollary of \cite{Nishio85}[Th.1, Corollary 3].
\begin{proposition}\label{prop1} Let $a,b \in \mathbb{R}\setminus \{0\}$,
 $P_1$, $P_2$ be orthoprojections on $H$. Then for every $x\notin
 \{0,a,b,a+b\}$ the following implications hold
 \[x\in \sigma(aP_1+bP_2)\Longleftrightarrow (a+b-x)\in
 \sigma(aP_1+bP_2).\]
 Beside this, both $x$ and $a+b-x$ have the same multiplicity as
 eigenvalues of $aP_1+bP_2$ or both are the approximate points of
 $\sigma(aP_1+bP_2)$.
  \end{proposition}
A simple application of Proposition \ref{prop1} leads to operator
inequalities for  linear combinations of two orthoprojections.
  \begin{cor}\label{cor2}
  Suppose there exists $c\in \mathbb{R}$ such that
  $aP_1+bP_2\le c I$. Then $(a+b-c)P_{\tilde{H}}\le aP_1+bP_2$,
  where $P_{\tilde{H}}$ is the orthogonal projection on the
  subspace $\tilde{H}=\overline{\Im P_1+\Im P_2}$.
  \end{cor}
By induction the statement of Corollary \ref{cor2} can be simply
expanded to the following proposition on operator inequalities
(see more discussion in \cite{RaYu10,Feshchenko12}).
\begin{proposition}\label{prop7} Let $\alpha_i>0$, $P_i$ be an
orthogonal projection in a Hilbert space $H$ for every $i=1,\dots
k$, $k\ge 2$. Suppose there exist $c\in \mathbb{R}$ such that the
following inequality holds: $\alpha_1 P_1+\dots +\alpha_k P_k\le c
I$.    Then \[(\alpha_1+\dots +\alpha_k-(k-1)c)P_{\tilde{H}}\le
\alpha_1 P_1+\dots +\alpha_k P_k,\]
  where $P_{\tilde{H}}$ is the orthogonal projection on the
  subspace $\tilde{H}=\overline{\Im P_1+\dots +\Im P_k}$.
\end{proposition}
For self-adjoint $A$, we shall denote by $\sigma_{ess}(A)$ the
Weyl's essential spectrum of $A$, that is the set
 \[\bigcap\limits_{K \in \mathcal{K}
 } \sigma(A+K)\]
where $\mathcal{K}$ is the set of all Hermitian compact operators.
 It consists of all limit points of $\sigma(A)$ and all
eigenvalues of $A$ of infinite multiplicity.
\begin{cor}
  In the setting of Proposition \ref{prop1} the following implications hold
 \[x\in \sigma_{ess}(aP_1+bP_2)\Longleftrightarrow (a+b-x)\in
 \sigma_{ess}(aP_1+bP_2).\]
  \end{cor}
 For more general results on combinations of
orthoprojections we refer the reader to nice surveys
\cite{Spitkovsky10,Wu94}.
\section{Main theorem}
In this section we consider direct sums of Hilbert spaces and
bounded operators on them. Let $H=V\oplus V$, where $V$ is a
separable Hilbert space. We denote by $L(H)$ the algebra of all
bounded operators on $H$. Every operator from $L(H)$ can be viewed
as $2\times 2$ block matrix. It is easy to show that for Hermitian
operators $T_1$, $T_2\in L(V)$ with $\sigma(T_1)\in [0,1]$ and
$\sigma(T_2)\in [0,2]$, the operator $\diag{T_1,-T_1}$ is a
difference of two orthoprojections and the operator
$\diag{T_2,2I_{V}-T_2})$ is a sum of two orthoprojections
\cite{Halmos2s}. For example defining orthoprojection $Q_1$ and
$Q_2$ by the formulas
\[
Q_1=\begin{pmatrix}(I_{V}+T_1)/2&\frac{1}{2}\sqrt{I_{V}-T_1^2}  \\
\frac{1}{2}\sqrt{I_{V}-T_1^2}&(I_{V}-T_1)/2
\end{pmatrix}, \quad Q_2=\begin{pmatrix}(I_{V}-T_1)/2&
\frac{1}{2}\sqrt{I_{V}-T_1^2}  \\
\frac{1}{2}\sqrt{I_{V}-T_1^2}&(I_{V}+T_1)/2
\end{pmatrix}, \]
we have $Q_1-Q_2=\diag{T_1,-T_1}$. Using this results and two
facts on self-commutators and linear combination of two
orthoprojection, we can prove the following theorem.
\begin{thm}\label{th1} Every bounded self-adjoint operator on $H$
is a linear combination of $4$ orthoprojections.
 \end{thm}
{\bf Proof}. 
Let $A$ be a self-adjoint operator. There
exist two subspaces $H_1$ and $H_2$ of $H$ such that $H=H_1\oplus
H_2$, $A H_1 \subset H_1$ and both $H_1$ and $H_2$ are isomorphic
to $H$. Without lost of generality we suppose that $H_2=H_1$.
According to this decomposition the operator $A$ has the block
diagonal form: $A=\diag{A_1,A_2}$, where $A_1$ and $A_2$ are
self-adjoint operators in $H_1$. Let $\lambda \in
\sigma_{ess}(A_1+A_2)$. Every Hermitian operator with $0$ in the
convex hull of its essential spectrum is a self-commutator
\cite{Radjavi66}. Since $0\in \sigma_{ess}(A_1+A_2-\lambda I)$,
there exits an operator $X$, such that $A_1+A_2-\lambda
I=[X^*,X]=X^*X-XX^*$. We can suppose that $X$ is invertible and
$X^*X> (|\lambda|+1)I$, because of the invariance property of the
commutator: $[X^*,X]=[X^*+tI,X+tI]$, $t\in \mathbb{C}$. Note that
$X^*X$ and $XX^*$ are unitary equivalent, so $\diag{X^*X, \lambda
I -XX^*}$ is a linear combination of $2$ orthoprojections
\cite{Nishio85}, say $aP_1-bP_2$ with $a,b> \|X^*X\|$ and $a-b
=\lambda$. Whence, we have
\begin{equation}\label{eq1} A-aP_1-bP_2=\diag{T,-T},
\end{equation}
with $T=A_1-X^*X$. Beside this, the operator $\diag{T,-T}$ is a
difference $P_3- P_4$ of two orthoprojections multiplied by
$\|T\|$. So for the number $c=\|T\|$,
\begin{equation}\label{eq3}
A=aP_1-bP_2+c P_3-c P_2, \end{equation}
 as required.\qed
\begin{cor}\label{cor4}
Every bounded operator on a Hilbert space is a linear combination
of $8$ orthoprojections.
    \end{cor}
 We call the orthoprojection $P$ \emph{proper}
if $\dim \Im P=\dim \Im (I-P)=\infty$.
\begin{cor}\label{cor3}
Every bounded self-adjoint operator $A$ on a Hilbert space is a
real linear combination of $5$ proper orthoprojections.
    \end{cor}
{\bf Proof}. 
 Suppose (\ref{eq3}) holds and some of $P_i$
are not proper. By construction, $cP_3-cP_4=\diag{T,-T}$. Whence
if $T$ is of infinite rank, then both $P_3$ and $P_4$ are proper.
For $T$ being of finite rank, we have $P_3$ and $P_4$ have a
common eigenspace $V_{34}$, $\dim V_{34}=\infty$ of the same
eigenvalue $\alpha \in \{0, 1\}$. Putting $P_{34}$ to be a proper
orthoprojection on a subspace of $V_{34}$ and
$\tilde{P}_i=P_i+(-1)^{\alpha}P_{34}$, $i=3,\ 4$, we see that
$\tilde{P}_3$ and $\tilde{P}_4$ are proper orthoprojections and
$\tilde{P}_3-\tilde{P}_4=P_3-P_4$.

Suppose now $P_1$ is proper and $P_2$ is not proper. If $\rank P_2
<\infty$ then it is a difference of two commuting orthoprojections
and if $\rank (I-P_2) <\infty$ then $P_2$ is a sum of two
commuting orthoprojections. Whence $A$ is a linear combination of
$5$ proper orthoprojections.

In view of symmetry it remains to consider only the case: both
$P_1$ and $P_2$ are not proper. Then they should have a common
eigenspace $V_{12}$, $\dim V_{12}=\infty$, such that $P_i
h=\alpha_i h$ for every $h \in V_{12}$, $i=1,\ 2$ and some
$\alpha_1$, $\alpha_2\in \{0,1\}$. Putting $P_{12}$ to be a proper
orthoprojection on a subspace of $V_{12}$ and
$\tilde{P}_i=P_i+(-1)^{\alpha_i}P_{34}$, $i=3,\ 4$, we see that
$\tilde{P}_1$ and $\tilde{P}_2$ are proper orthoprojections and
$A-(a\tilde{P}_1-b\tilde{P}_2+cP_3-cP_4)=
((-1)^{\alpha_2}b-(-1)^{\alpha_1}a)P_{12}$. Hence $A$ is a linear
combination of $5$ proper orthoprojections.\qed
\begin{rem}\label{rem1}
In the proof of Theorem \ref{th1}, we can take $\lambda =0$ in
case the number $0$ is in the convex hull of $\sigma(A_1+A_2)$ and
then put $a=b\in \mathbb{N}$ and set $c$ being integer part of
$||T||+1$. Thus $A$ is an integral linear combination of $4$
orthoprojections for this particular case. In general situation,
when  $\lambda\in\sigma_{ess}(A_1+A_2)$ and $\lambda\ne 0$, we can
take any proper orthoprojection from $L(H_1)$ and any integer
number $d$ satisfied the conditions $|d|
> 2\|A\|$ and $\lambda d>0$. Then we obtain that convex hull of
$\sigma_{ess}(A_1+A_2-dP_5)$ contains zero. So the operator $A-d \
\diag{P_5,0_{H_1}}$ is an integral linear combination of $4$
orthoprojections and $A$ is an integral linear combination of $5$
orthoprojections.
\end{rem}
\section{Counterexamples}\label{sec4}
Not every operator of the form $I+K$, where $K$ is an
infinite-rank compact operator is a linear combination of $3$
orthoprojections. To prove this we have to show that for $P_1$
being orthoprojection, the difference $I+K-\alpha_1 P_1$ is not a
linear combination of two orthoprojection. Note that for proper
$P_1$ with $\alpha_1\ne 1$, the spectrum
$\sigma_{ess}(I+K-\alpha_1 P_1)$ has exactly two points and if
$P_1$ is not proper, then the corresponding essential spectrum has
only one point. So we start with properties of linear combinations
of two orthoprojections whose essential spectra contain at most
two points.

\emph{(i)} Let $a\le b\in \mathbb{R}$,  $P_1$, $P_2$ be
orthoprojections and suppose that $\sigma_{ess}(aP_1+b
P_2)=\{x\}$. Then $x\in \{0, a,b,a+b\}$. Assume that this is not
true, i.e. $x\notin \{0, a,b,a+b\}$. Applying Corollary
\ref{cor2}, we obtain $a+b-x\in\{x\}\Longrightarrow x=(a+b)/2$. In
view of (\ref{eq2}) this can be true only if $a=b$ or $a+b=0$. So
$x=a$ or $x=0$. A contradiction.

 We note that for a sequence of
different numbers $x_1,x_2,x_3,\dots$ with every $x_i\in
\sigma(aP_1+b P_2)$ and $\lim_{i\to\infty} x_i=x$, we have here by
Proposition \ref{prop1} that all except may be two elements of the
sequence $a+b-x_1, a+b-x_2,\dots $ must be in $\sigma(aP_1+b P_2)$
and so the spectrum contains infinite number of points less than
$x$ and infinite number of points greater than $x$.

\emph{(ii)} Let  now $\sigma_{ess}(aP_1+b P_2)=\{x, y\}$, where
$0<x<y$ and $0<a\le b$. Then $x+y = a+b$ or $x\in \{a,b\}$.
Indeed, suppose the inverse is true, i.e. $x+y \ne a+b$ and
$x\notin \{a,b\}$. Since $x<y$ and $y\in \sigma(aP_1+b P_2)$, then
$x<a+b$. By Corollary \ref{cor2}, $a+b-x\in \sigma_{ess}(aP_1+b
P_2)$ and so $a+b-x=x$ or $a+b-x=y$. The last equality is not
valid by assumption, hence $x=(a+b)/2=a$ as in the previous
paragraph. A contradiction.

 We note that $a+b\ne x+y$ yields also to the
equality $y=a+b$ and the number $a+b$ has to be is an isolated
point of $\sigma(aP_1+b P_2)$. Indeed, assuming the existence of
different $x_i\in \sigma(aP_1+b P_2)$, $\lim_{i\to \infty}x_i=
a+b$, we see $\lim_{i\to \infty}(a+b-x_i)=0$, whence by
Proposition \ref{prop1}, $0\in \sigma_{ess}(aP_1+b P_2)$, which is
not true by the initial assumption.

We shall use the following result on rank-one perturbation of a
Hermitian compact operator \cite{Hochstadt73}.
\begin{proposition}\label{prop2}
Let $K=\diag{\mu_1,\mu_2,\mu_3,\dots}$, where
$\mu_1>\mu_2>\mu_3>\dots$, $\mu_n\to 0$, $n\to \infty$. For every
rank one orthogonal projection $P$ and   every $t>0$, the set of
eigenvalues $\gamma_1\ge \gamma_2\ge\gamma_3\ge \dots $ of $K+tP$
satisfy the interlace relation $\gamma_1\ge \mu_1 \ge \gamma_2\ge
\mu_2 \ge  \gamma_3\ge \dots$
\end{proposition}
We note that for $t<0$ the interlace property in the proposition
is also true but in the inverse order: $\dots,\gamma_3\le
\mu_3\le\gamma_2\le \mu_2 \le \gamma_1\le \mu_1 $, see also more
discussion in \cite{Vasudeva76}.
\begin{proposition}\label{prop3}
Let $K$ be non-negative compact operator of infinite rank. Then
$I-K$ is not a linear combination of three orthoprojections.
\end{proposition}
{\bf Proof}. 
 The main goal in the proof is to find a
relation between coefficients in a decomposition of $I-K$ into a
linear combination of $3$ orthoprojections if such a decomposition
exists. And then prove that with such a relation the decomposition
do not exist. Let us split a proof into several parts.

1. We assume $I-K=\beta_1Q_1+\beta_2Q_2+\beta_3Q_1$, where
$\beta_1,\beta_2,\beta_3\in \mathbb{R}$ and $Q_1$, $Q_2$ and $Q_3$
are orthoprojections. If some of the coefficients, say $\beta_1$
and $\beta_2$, are negative, then we have \begin{align}
 I-K=&
\beta_1Q_1+\beta_2Q_2+\beta_3Q_1 \Longleftrightarrow \notag \\
(1-\beta_1-\beta_2)I-K=&(-\beta_1)(I-Q_1)+(-\beta_2)(I-Q_2)+
\beta_3Q_1.\notag
\end{align}
 So
$(1-\beta_1-\beta_2)I-K$ is a linear combination of $3$
orthoprojections with positive coefficients. Putting
$\tilde{K}=K/(1-\beta_1-\beta_2)$,  we conclude that there exists
decomposition
\[I-\tilde{K}=\alpha_1P_1+\alpha_2P_2+\alpha_3P_3,\]
where $P_1$, $P_2$, $P_3$ are orthoprojections, $0\le \alpha_1\le
\alpha_2\le \alpha_3\le 1$ and $\tilde{K}\ge 0$ is a compact
self-adjoint operator of infinite rank. We denote
$B_{ij}=\alpha_iP_i+\alpha_jP_j$ for $i\ne j$. The coefficient
$\alpha_1>0$ or otherwise $\alpha_1=0$ and $I-\tilde{K}=B_{23}$ is
a linear combination of two orthoprojection. Immediately we have
$\sigma_{ess}(B_{23})=\{1\}$. This is the case (i) above with
$x=1$. But $1$ is an approximate point of $\sigma(B_{23})$ and all
but one points of $\sigma(B_{23})$ is less then $1$ by definition.
So $I-\tilde{K}$ can not be a linear combination of two
orthoprojections. Whence $\alpha_1>0$.

The coefficient $\alpha_3< 1$. Indeed, if $\alpha_3=1$, then
$B_{12}=(I-P_3)-\tilde{K}\ge 0$. Whence $\ker (I-P_3)\subset \ker
\tilde{K}$. Considering restrictions of $P_1$ and $P_2$ onto
$\tilde{H}=Im (I-P_3)$, we see that restriction
$B_{12}\left|_{\tilde{H}}\right.$ is a linear combination of two
orthoprojections and at the same time
$B_{12}\left|_{\tilde{H}}\right.
=I_{\tilde{H}}-\tilde{K}\left|_{\tilde{H}}\right.$ with
$\tilde{K}\left|_{\tilde{H}}\right.$ of infinite rank.
 So we come to a contradiction as in the case $\alpha_1=0$ above.

2. Now we can prove that every $P_i$ is proper, $i=1,2,3$. For
example, if $P_1$ is of finite rank, then $-\tilde{K}- \alpha_1
P_1$ is non-positive compact operator of infinite rank. So
$B_{23}=I-(\tilde{K}+ \alpha_1 P_1)$ can not be a linear
combination of two orthoprojections (see the explanation of the
case $\alpha_1=0$).

If $I-P_1$ is of finite rank, then
$B_{23}=(1-\alpha_1)I+\alpha_1(I-P_1)-\tilde{K}$. Hence
$\sigma_{ess}(B_{23}) =\{1-\alpha_1\}$ and the compact operator
$\alpha_1(I-P)-\tilde{K}$ is finite dimensional perturbation of
non-positive operator $-\tilde{K}$. So it has finite number of
positive eigenvalues and infinite number of negative eigenvalues.
This means $1-\alpha_1$ is an approximate point of
$\sigma(B_{23})$ by infinite different numbers less than
$1-\alpha_1$ and there are only finite numbers from
$\sigma(B_{23})$ greater than $1-\alpha_1$. Hence $B_{23}$ can not
be a linear combination of $2$ orthoprojections or $P_{1}$ has to
be proper. In view of symmetry, $P_i$ is proper for every
$i=1,2,3$.

3. Let us consider $\sigma_{ess}(B_{23})$ more closely. The
operator $I-\alpha_1 P_1$ has two points of infinite multiplicity
in its spectrum, $1$ and $1-\alpha_1$, so
$\sigma_{ess}(B_{23})=\{1,1-\alpha_1\}$. There are only two
possible values for the sum $\alpha_2+\alpha_3$:
$\alpha_2+\alpha_3=1$ and $\alpha_2+\alpha_3=2-\alpha_1$. We
consider these cases separately.

Let $\alpha_2+\alpha_3=1$, i.e. $1$ is an isolated point of
$\sigma(B_{23})$. We define $H_{23}=(\Im P_2\cap \Im
P_3)^{\perp}$. For every $h\in H_{23}^{\perp}$,
$P_2h=P_3h=h\Longrightarrow (\alpha_2P_2+\alpha_3 P_3)h=h$. Also
$((I-\tilde{K})x,x)\le \|x\|^2$ for every $x\in H$. On the other
hand,
$((I-\tilde{K})h,h)=(\alpha_1P_1+\alpha_2P_2+\alpha_3P_3h,h)=\alpha_1(P_1h,h)+
\|h\|^2$. This yields $P_1h=0$. So $ H_{23}^{\perp}\in \Ker P_1$
and $ H_{23}^{\perp}\in \Ker \tilde{K}$. Therefore we can restrict
every operator $P_1$, $P_2$, $P_3$ and $\tilde{K}$ to $H_{23}$,
obtaining the decomposition
$\hat{I}-\hat{K}=\alpha_1\hat{P}_1+\alpha_2\hat{P}_2+\alpha_3\hat{P}_3$
on Hilbert space $H_{23}$. In this decomposition $\hat{K}$ is
obviously of infinite rank but the orthoprojections $\hat{P}_1$,
$\hat{P}_2$ and $\hat{P}_3$ might not be proper. We repeat the
same argument from the part $2$ of the proof and so we can assume
without lost of generality that every $\hat{P}_{i}$ is proper,
$i=1,2,3.$ Again $\alpha_2\hat{P}_2+\alpha_3\hat{P}_3=
(\hat{I}-\alpha_1\hat{P}_1)-\hat{K}$, hence $1$ must be in
$\sigma_{ess}(\alpha_2\hat{P}_2+\alpha_3\hat{P}_3)$ but by
construction $1\notin
\sigma(\alpha_2\hat{P}_2+\alpha_3\hat{P}_3)$. Therefore
$\alpha_2+\alpha_3\ne 1$.

4. The remaining case is $\alpha_1+\alpha_2+\alpha_3=2$. We remind
that $\alpha_j<1$ and $\alpha_i\le \alpha_{i+1}$, hence
$\alpha_2+\alpha_3\ge 4/3$. In view of Corollary \ref{cor2} from
the inequality $B_{23}\le I$, we have $(\alpha_2+\alpha_3-1)P_{\Im
P_2+\Im P_3}\le B_{23}$, that is $x\in
\sigma(B_{23})\Longrightarrow x=0$ or $x\ge \alpha_2+\alpha_3-1$.
So $0$ is an isolated point of $B_{23}$.

 Let
$\tilde{K}$ be a diagonal operator $
\diag{\gamma_1,\gamma_2,\gamma_3,\dots}$ in some orthogonal base
of $H$ and $\gamma_i\ge \gamma_{i+1}$, $i\ge 1$. There exist $k\in
\mathbb{N}$ such $\gamma_{k+1}<1-\alpha_1$. Since the operator
$I-\alpha_1P_1\ge (1-\alpha_1)I$, the operator
$I-\alpha_1P_1-\diag{0_k,\gamma_{k+1},\gamma_{k+2},\dots}$ is
invertible. Hence  dimension $n$ of the kernel of $B_{23}$ is at
most $k$. Let $V_{23}=(\Ker B_{23})^{\perp}$. We define
 $\check{P}_1$ as a  rank $n$ orthoprojection onto a subspace of
 $\Im P_1$ with the following property: $(P_1-\check{P}_1)h=0$ for
 every $h \in V_{23}^{\perp}$. Putting
 $\check{K}=\tilde{K}+\alpha_1\check{P}_1$, we have
 \[I-\check{K}=\alpha_1(P_1-\check{P}_1)+\alpha_2P_2+\alpha_3P_3.\]
By construction, every operator $P_1-\check{P}_1$, $P_2$ and $
P_3$ maps $V_{23}^{\perp}$ into zero vector. Hence
$V_{23}^{\perp}$ is invariant under the act of these operators and
so is $\check{K}$. As in the previous part we consider the
restriction of the operators to the subspace $V_{23}$ marking
corresponding operators with breve:
$\breve{I}-\breve{K}=\alpha_1\breve{P}_1+\alpha_2\breve{P}_2+\alpha_3\breve{P}_3.$
Since here $\alpha_2\breve{P}_2+\alpha_3\breve{P}_3$ is
invertible, we have $\breve{I}-\breve{K}-\alpha_1\breve{P}_1$ is
invertible. On the other hand,
\begin{equation}\label{eq54}\alpha_2\breve{P}_2+\alpha_3\breve{P}_3\ge
(\alpha_2+\alpha_3-1)\breve{I}=(1-\alpha_1)\breve{I}.
\end{equation}
 Also, according to the decomposition of $V_{23}$
into a direct  sum,
\begin{equation}\label{eq55}
V_{23}=H_1\oplus H_2, \quad H_2 = \Im \breve{P}_1,
\end{equation}
the operator $\breve{I}-\alpha_1\breve{P}_1$ has the diagonal form
$\diag{I_{H_1},(1-\alpha_1)I_{H_2}}$.  The inequality (\ref{eq54})
implies $\breve{I}-\breve{K}-\alpha_1\breve{P}_1\ge
(1-\alpha_1)\breve{I}$, whence $\breve{K}H_2=0$ and so
$\breve{K}=\diag{K',0_{H_2}}$ subject to the decomposition
(\ref{eq55}). As a corollary we have that $1$ is an approximated
point of $\sigma(\alpha_2\breve{P}_2+\alpha_3\breve{P}_3)$ and
this point is greater than $\alpha_3$ and less than
$\alpha_2+\alpha_3$ and at the same time
$\alpha_2+\alpha_3-1=1-\alpha_1$ is an isolated point of
$\sigma(\alpha_2\breve{P}_2+\alpha_3\breve{P}_3)$. Therefore
$\breve{I}-\breve{K}-\alpha_1\breve{P}_1$ can not be a linear
combination of two orthoprojections and this complete the proof of
the part.
  \qed
\begin{cor}\label{cor45}
Let $K$ be non-negative compact operator of infinite rank. Then
$I+K$ is not a linear combination of three orthoprojections.
\end{cor}
{\bf Proof}. 
 Suppose $I+K=\beta_1Q_1+\beta_2 Q_2+\beta_3
Q_3$, where for every $i=1,2,3$, $Q_i$ is an orthoprojection.
Replacing $Q_i$ with $I-Q_i$ when $\beta_i<0$, we can find new
decomposition $I+c K=\alpha_1 P_1+ \alpha_2 P_2+ \alpha_3 P_3$
with $\alpha_i\ge 0$ and some positive $c$. The equivalent
decomposition is
\begin{equation}\label{eq411}
(\sum\limits_{i=1}^3 \alpha_i -1)I-c K= \sum\limits_{i=1}^3
\alpha_i (I-P_i).\end{equation} Since
$\alpha_1+\alpha_2+\alpha_3\ge \|I+c K\|>1$, the decomposition
(\ref{eq411}) states that scalar operator minus compact operator
is a linear combination of $3$ orthoprojections, which contradicts
Proposition \ref{prop3}. This completes the proof.
  \qed

We denote by $\boldsymbol{i}$ the imaginary unit $\sqrt{-1}$ in
$\mathbb{C}$, by $Re(x)$ and $Im(x)$ corresponding real and
imaginary parts of a complex number $x$.
\begin{cor}\label{cor5} The operator of the form $I-K-\boldsymbol{i}
K$ is not a complex linear combination of $4$ orthoprojections.
    \end{cor}
{\bf Proof}. 
 Suppose the inverse and $I-K-\boldsymbol{i}
K=c_1P_1+c_2P_2+c_3P_3+c_4P_4$. Then
\begin{equation}\label{eq10}
I-K=Re(c_1)P_1+Re(c_2)P_2+Re(c_3)P_3+Re(c_4)P_4 \end{equation}
 and
\begin{equation}\label{eq11}
-K=Im(c_1)P_1+Im(c_2)P_2+Im(c_3)P_3+Im(c_4)P_4. \end{equation}
 By Proposition
\ref{prop3}, every $Re(c_i)\ne 0$, $i=1,\dots,4$. Also
\[-K-\sum\limits_{Im(c_j)<0} Im(c_j)I =
\sum\limits_{1}^{4}\left| Im(c_i)\right| Q_i, \quad Q_i=
  \begin{cases}
    P_i & \text{ if } Im(c_i)\ge 0, \\
   I-P_i& \text{otherwise}.
  \end{cases}
\]
Thus for every $i=1,\dots,4$, $Im(c_i)\ne 0$ either.  There exists
$k$, such that $Re(c_k)\ne Im(c_k)$, because $I-K\ne -K$.
Evaluating $I-K-Re(c_k)/ Im(c_k)K$ by (\ref{eq10}) and
(\ref{eq11}), we have that this operator is a real linear
combination of orthoprojections $P_i$, $i=1,\dots,4,$ $i\ne k$,
which is not true by Proposition \ref{prop3} or by Corollary
\ref{cor45}. Therefore $I-K-\boldsymbol{i} K$ is not a linear
combination of $4$ orthoprojections.
 \qed

Now we turn our attention to finite matrices. Impossibility to
decompose $I-K$ into a linear combination of $3$ orthoprojections
from Proposition \ref{prop3} suggests the form of a matrix for
which such a decomposition does not exist either. Before we start
we recall that the interlace property from Proposition 2 is also
true for Hermitian matrices. Following \cite{HJ85} we denote by
$\lambda_k(A)$ the $k$-th smallest eigenvalue of the Hermitian
matrix $A$ counting multiplicity. We shall use well known Weyl's
theorem on rank $k$ perturbation of spectrum of a Hermitian matrix
(see \cite[Theorem 4.3.6]{HJ85})
\begin{thm}\label{th5} Let $A$ and $B$ be Hermitian $n\times m$
matrices and $\rank B \le k$. Then $\lambda_j(A+B)\le
\lambda_{j+k}(A)\le \lambda_{j+2k}(A+B)$, $j=1,\dots, n-2k$.
\end{thm}
Also we shall frequently use the monotonicity property for
eigenvalues of Hermitian $n\times n$ matrices: $A\le
B\Longrightarrow \lambda_i(A)\le \lambda_i(B)$ for every
$i=1,\dots,n$.
\begin{proposition}\label{prop6}
Let $A=\diag{\mu_1,\dots,\mu_4,\gamma_1I_{18},\dots,\gamma_4
I_{18}}$, where $\mu_i =(1-  10^{-10i}\theta)$, $\gamma_i=(1+
10^{100(i-5)}\theta)$, $i=1,\dots,4$, $0<\theta\le 1$. Then $A$ is
not a real linear combination of three orthoprojections.
\end{proposition}
{\bf Proof}. 
The proof is by contradiction. The matrix
$A$ is not a linear combination of two orthoprojection by
Proposition \ref{prop1}.  So suppose $A$ is a linear combination
of $3$ orthprojections, say $
A=\alpha_1P_1+\alpha_2P_2+\alpha_3P_3$. Let $n$ be the size of
$A$, $n=76$. If one of the coefficients $\alpha_1$, $\alpha_2$,
$\alpha_3$ is negative, then using the procedure from the first
part of the proof of Proposition \ref{prop3}, we come to a new
matrix $A_1$ and the decomposition with positive coefficients,
\[A_1=\frac{A+cI_n}{1+c}=\sum\limits_{i=1}^3 \frac{|\alpha_i|}{1+c}
\tilde{P}_i,\] where
$c=(|\alpha_1|+|\alpha_2|+|\alpha_3|-\alpha_1-\alpha_2-\alpha_3)/2$.
Note that in this case eigenvalues of $A_1$ can be calculated by
formulas for eigenvalues of $A$ in the formulation of Proposition
\ref{prop6} but with smaller value of parameter $\theta$, which
should be equal $\theta/(1+c)$.  We will not specify the parameter
$\theta$ and so, without lost of generality, we may assume that
every $\alpha_i$ is positive, $i=1,2,3$ and $\alpha_1\le
\alpha_2\le \alpha_3\le \gamma_4$. We denote $\epsilon=1-\mu_1$,
$\delta=\gamma_1-1$, $B_{ij}=\alpha_iP_i+\alpha_jP_j$ and
$H_{ij}:=\Im B_{ij}$.

Let us show that $9\le \rank P_i \le n-9 $ for every $i=1,2,3$.
It's enough to establish this for $P_1$. At first, suppose
$k_1=\rank P_1<9$. Since $B_{23}$ is a rank $k_1$ perturbation of
$A$, we have that every eigenvalue of $B_{23}$ which does not
coincide with $\gamma_1,\dots,\gamma_4$ has multiplicity at most
$k_1+1\le 9$ and every $\gamma_i$ has multiplicity at least
$18-k_1\ge 10$. Since $B_{23}$ is a linear combination of $2$
orthoprojections, then by Proposition \ref{prop1}, $\gamma_1$
coincides with one of the numbers $\alpha_2$, $\alpha_3$ and
$\alpha_2+\alpha_3$ or the number $\alpha_2+\alpha_3-\gamma_1$ is
an eigenvalue of $B_{23}$ of the same multiplicity as $\gamma_1$
and so it coincides with $\gamma_j$ for some $j=1,\dots, 4$. The
same is true for $\gamma_2$, $\gamma_3$ and $\gamma_4$. Therefore,
there exist $i_1$, $i_2$, $i_3$, $i_4$ such that $\{i_1,i_2\}\ne
\{i_3,i_4\}$ and
\begin{equation}\label{eq12}
\gamma_{i_1}+\gamma_{i_2}=\alpha_2+\alpha_3=\gamma_{i_3}+\gamma_{i_4}.
\end{equation}
 In view of definition of $\gamma_i$, the property (\ref{eq12})
 does not hold for any different sets $\{i_1,i_2\}$ and
 $\{i_3,i_4\}$. So $\rank P_1\ge 9$.

Suppose now that $\rank P_1>n-9$. 
Then $B_{23}=(A-I)+(I-P_1)$ and hence $B_{23}$ is a rank $n-k_1$
perturbation of $A-I$ with $n-k_1<9$. The same reason as above
shows $\gamma_{i_1}+\gamma_{i_2}=\gamma_{i_3}+\gamma_{i_4}$, which
is not true by definition. So $\rank P_1\le n-9$. In view of
symmetry $9\le \rank P_i \le n-9 $ for every $i=1,2,3$.

Thus $9\le k_1\le n-9$. By Theorem \ref{th5} for the matrix
$B_{23}$ and $\alpha_1P_1$, we have $\gamma_{1}= \lambda_{5}(A)\le
\lambda_{k_1+5}(B_{23})$. We fix some $p\in \mathbb{N}$ with
$\lambda_{p}(B_{23})\ge \gamma_1$ and define
$x^*=\lambda_{p}(B_{23})$. Since $B_{23}$ is a linear combination
of $2$ orthoprojections, there exist only four possible cases for
$x^*$:
\begin{align}
1)\ & x^*=\alpha_2, \quad 2)\ x^*=\alpha_3, \quad 3)\
x^*=\alpha_2+\alpha_3, \notag\\ 4)\ & x^*\notin
\{\alpha_2,\alpha_3,\alpha_2+\alpha_3\}, \quad \alpha_2+\alpha_3-
x^*\in \sigma(B_{23}).\notag
\end{align}
We consider all cases separately.

\emph{Case 1)} \ $x^*=\alpha_2$. We use only the fact that
$\alpha_2\ge \gamma_1$.   Due to ordering, $\alpha_3\ge
\alpha_2\ge \gamma_1$. Applying Corollary \ref{cor2} to
inequalities $B_{23}\le A\le \gamma_4 I$, we get
\begin{equation}\label{eq13}
(\alpha_2+\alpha_3-\gamma_{4})P_{H_{23}}\le \alpha_2
P_2+\alpha_3P_3.
\end{equation}
So all nonzero eigenvalues of $B_{23}$ is greater or equal to
$\alpha_2+\alpha_3-\gamma_{4}\ge 2\gamma_1-\gamma_4>\mu_4$. Since
$B_{23}\le A$, we have by monotonicity principle,
$\lambda_4(B_{23})\le \mu_4$. In view of (\ref{eq13}), we obtain
that $\lambda_1(B_{23})=\dots =\lambda_4(B_{23})=0$. On the other
hand, $\alpha_1P_1+B_{23}=A$, so $\alpha_1\ge \mu_4$. Applying now
Proposition \ref{prop7} to the linear combination
$\alpha_1P_1+\alpha_2P_2+\alpha_3P_3$, we conclude
\begin{equation}\label{eq14}
A\ge (\alpha_1+\alpha_2+\alpha_3-2\gamma_4)I\ge
(2\gamma_1+\mu_4-2\gamma_4)I> \mu_3 I, \end{equation}
 that is
$\mu_3 \notin \sigma(A)$ which is not true by definition of $A$.
Therefore $\alpha_2<\gamma_1$ and so $x^*\ne \alpha_2$.

\emph{Case 2)} \ $x^*=\alpha_3$. Here we consider
$B_{12}=A-\alpha_3P_3$. Since $\rank P_3\le n-9$, then by Theorem
\ref{th5} we have $\lambda_{n}(B_{12})\ge \lambda_9(A)= \gamma_1$.
So $\alpha_1+\alpha_2\ge \gamma_1$ and $\alpha_2\ge \gamma_1/2$.
Putting $\triangle=\gamma_4-\gamma_1$ and substituting $\gamma_1$
for $\alpha_3$ in (\ref{eq13}), we obtain
\begin{equation}\label{eq15}
B_{23}\ge (\alpha_2-\triangle)P_{H_{23}}.
\end{equation}
 If in
addition $P_{H_{23}}\ne I$, then $B_{23}$ is singular and so
$B_{23}+\alpha_1 I$ has an eigenvalue $\alpha_1$. With the
property $A\le B_{23}+\alpha_1 I$ this implies by monotonicity
property that $\alpha_1\ge \mu_1$ and so automatically
$\alpha_2\ge \mu_1$. In view of (\ref{eq15}) and
$B_{23}+\alpha_1P_1\le A$, we conclude $\Im P_1\cap
H_{23}=\emptyset$. So $\dim \ker (B_{23})=\rank P_1\ge 9$. Hence
$B_{23}+\alpha_1 I$ has the eigenvalue $\alpha_1$ of multiplicity
at least $9$ that is $\alpha_1=\lambda_{9}(B_{23}+\alpha_1 I)\ge
\lambda_{9}(A)=\gamma_1$. As a corollary we obtain
$\alpha_2\ge\gamma_1$ and this is case 1) which  was considered
above.

Thus, $P_{H_{23}}= I$. It follows immediately from (\ref{eq15})
that $(\alpha_2-\triangle)I+\alpha_1P_1\le A$. So
$\alpha_1+\alpha_2-\triangle\le \gamma_4$. We mentioned in
Preliminaries that the spectrum of a linear combination
$\alpha_2P_2 +\alpha_3P_3$ of orthoprojections lies in the union
of two segments: $[0,\alpha_2]\cup [\alpha_3,\alpha_2+\alpha_3]$.
Hence by (\ref{eq15}), we have $\sigma(B_{23})\cap
[0,\alpha_2]\subset [\alpha_2-\triangle,\alpha_2]$. This yields
from Proposition \ref{prop1}, that $\sigma(B_{23})\cap
(\alpha_2,\alpha_2+\alpha_3]\subset
[\alpha_3,\alpha_3+\triangle]$. Hence $B_{23}=\alpha_3 (I-Q)+
\alpha_2 Q + B_{\triangle}$, where $Q$ is some orthogonal
projection and $B_{\triangle}$ is a Hermitian matrix with
$\|B_{\triangle}\| \le {\triangle}.$  We note that
\begin{align}\label{eq16}
\alpha_1P_1+(\alpha_3-\alpha_2)Q=\alpha_1P_1+\alpha_3 Q +
\alpha_2(I-Q)-\alpha_2I  \le \alpha_1P_1 \notag \\
+B_{23}-\alpha_2 I +\triangle I \le (\gamma_5-\alpha_2+\triangle)I
\le (\alpha_3-\alpha_2+2\triangle)I. \end{align}
 By Corollary
\ref{cor2} for linear combination of orthoprojections $P_1$ and
$Q$, we have
\begin{equation}\label{eq17}
(\alpha_1-2\triangle)P_{H_1}\le \alpha_1P_1+(\alpha_3-\alpha_2)Q,
\end{equation} where $H_1=\Im P_1+\Im Q$. As we showed above,
$\alpha_1+\alpha_2\ge \gamma_1$ and $ \gamma_1\le\alpha_3\le
\gamma_4$. Hence $\alpha_1\ge \gamma_1-\alpha_2\ge \alpha_3
-(\gamma_4-\gamma_1) -\alpha_2=\alpha_3-\alpha_2-\triangle$.
Combining (\ref{eq16}) and (\ref{eq17}), we obtain
\begin{equation}\label{eq18}
(\alpha_3-\alpha_2-3\triangle)P_{H_1}\le
\alpha_1P_1+(\alpha_3-\alpha_2)Q\le
(\alpha_3-\alpha_2+2\triangle)I,
\end{equation}
that is $\sigma(\alpha_1P_1+(\alpha_2-\alpha_3)Q)\in \{0,
[\alpha_3-\alpha_2-3\triangle,\alpha_3-\alpha_2+2\triangle]\}.$ By
construction, $A=\alpha_1P_1+(\alpha_3-\alpha_2)Q+ \alpha_2
I+B_{\triangle}$. So the spectrum of $A$ must be in $4\triangle$
neighborhoods of the following three points $0$, $\alpha_2$ and
$\alpha_3$. By conditions of the Proposition, the eigenvalues
$\mu_i$ of $A$ satisfies the inequalities
$|\mu_i-\mu_j|>8\triangle $ for $i\ne j$. So $\sigma(A)$
 contains a point that is not from the mentioned neighborhoods and
 therefore $x^*\ne \alpha_3$.

\emph{Case 3)} \ $x^*=\alpha_2+\alpha_3$. We note that in this
case $\alpha_3\ge x^*/2$ and $\alpha_2\le x^* /2$. Hence
$\alpha_1\le x^*/2$. It follows then
\begin{equation}\label{eq19}
B_{23}=A-\alpha_1 P_1\ge \mu_1 I-\alpha_1 I\ge (\mu_1-x^*/2)I,
\end{equation}
i.e. $B_{23}$ is invertible. On the other hand, $B_{23}\le
\alpha_3 P_3+\alpha_2 I$ and the last matrix has at least five
pairwise orthogonal  eigenvectors with eigenvalue $\alpha_2$ since
$\rank P_3\le n-9$. This means $\lambda_5(B_{23})\le \alpha_2$.
Let $H_{\mu}$ be four dimensional subspace containing eigenvectors
of $A$ corresponding the eigenvalues $\mu_1$, $\mu_3$, $\mu_3$ and
$\mu_4$. Applying Courant-Fischer min-max theorem for the
eigenvalue $\lambda_5(B_{23})$, we get
\begin{align}\label{eq20}
\alpha_2& \ge  \lambda_5(B_{23}) \ge  \min\limits_{\|x\|=1, x\perp
H_{\mu} } (B_{23} v,v) \notag \\ &\ge  \min\limits_{\|v\|=1,
v\perp H_{\mu} } (A v,v) -\max\limits_{\|v\|=1, v\perp H_{\mu} }
(\alpha_1P_1 v,v)= \gamma_1-\alpha_1,
\end{align}
where $(v,w)$ means the inner  product of corresponding vectors.
So $\alpha_1+\alpha_2\ge \gamma_1$. Also $\gamma_4\ge
x^*=\alpha_2+\alpha_3\ge 2\alpha_2 \ge \alpha_2+\alpha_1$.
Combining these two inequalities, we obtain estimations on
$\alpha_1$, $\alpha_2$ and $\alpha_3$:
\begin{description}
\item[$\bullet$]
 $\gamma_1 /2\le \alpha_2 \le \gamma_5/2$,
\item[$\bullet$]
$\gamma_1/2\le \alpha_3 \le \gamma_5-\gamma_1/2=\gamma_5 /2
+\triangle/2$,
\item[$\bullet$] $\gamma_1 /2 -\triangle/2=
\gamma_1-\gamma_5/2\le \alpha_1 \le \gamma_5/2$.
\end{description}
Thus the numbers $\alpha_1$, $\alpha_2$ and $\alpha_3$ are in
$\triangle$ neighborhood of $1/2$. Using Proposition \ref{prop1}
and (\ref{eq19}) to $\alpha_2P_2+\alpha_3P_3$, we get
\begin{equation}\label{eq21}
\sigma(B_{23})\subset \{ [\mu_1-x^* /2, \alpha_2+\alpha_3- \mu_1+
x^* /2], \alpha_2+\alpha_3 \}\subset \{[0.4,0.6],
\alpha_2+\alpha_3 \}. \end{equation}
 Let $H_2=\Im P_2\cap \Im
P_3$, the eigensubspace of $B_{23}$ corresponding the eigenvalue
$\alpha_2+\alpha_3$. Since $\mu_1 I\le A\le \gamma_5 I$, then
$\mu_1 I -\alpha_1P_1\le A-\alpha_1P_1\le \gamma_4 I
-\alpha_1P_1$. Using monotonicity property, we conclude
that$A-\alpha_1P_1$ has at least $k_1$ eigenvalues which are not
greater than $\gamma_4-\alpha_1$ and  at least $n-k_1$ eigenvalues
which are greater or equal $\mu_1$. So
$\lambda_{k_1}(B_{23})\le\mu_1\le \lambda_{k_1+1}(B_{23})$. In
view of (\ref{eq21}), we have $\lambda_{k_1+1}(B_{23})=\dots=
\lambda_{n}(B_{23})$ and $\dim H_2=n-k_1$. We remind that $B_{23}$
is invertible. So $\rank P_2+\rank P_3=n+\dim H_2=2n-k_1$. Let us
estimate the trace of $A$:
\begin{equation}\label{eq211}tr A=\sum\limits_1^n \lambda_i(A)\le
(n-4)\gamma_4 +\mu_1+\mu_2+\mu_3+\mu_4 \le
n\gamma_4-(\gamma_4-\mu_1).
\end{equation}
Also, \begin{equation}\label{eq212}
tr A=\sum\limits_1^3 \tr
(\alpha_i P_i)=\sum\limits_1^3 \alpha_i \rank P_i
 \ge \alpha_1 \sum\limits_1^3  \rank P_i =2n \alpha_1\ge
 n(\gamma_1-\triangle).
   \end{equation}
Subtracting (\ref{eq212}) from (\ref{eq211}), we get
\[(n\gamma_5-(\gamma_5-\mu_1))-n(\gamma_1-\triangle)=(n+1)\triangle
- (\gamma_5-\mu_1)\ge 0.\]
 This inequality is not valid by
conditions of the Proposition. Therefore $x^*\ne
\alpha_2+\alpha_3$.

Before we start with case 4) we remark that according to the proof
of the cases 1)--3) the number $x^*$ can not belong to $\{
\alpha_2,\alpha_3, \alpha_2+\alpha_3\}$ for every $p$. So without
lost of generality we assume that every eigenvalue $\lambda$ of
$B_{23}$ does not belong $\{ \alpha_2,\alpha_3,
\alpha_2+\alpha_3\}$ as soon as it is greater or equal $\gamma_1$.
Also in case 1) we prove more strong statement that in all
possible decompositions of $A$ into a linear combination of $3$
orthoprojection with positive coefficients, the value of
$\alpha_2$ has to be less than $\gamma_1$.

 \emph{Case 4)} \ $x^*\notin\{
\alpha_2,\alpha_3, \alpha_2+\alpha_3\}$, $\alpha_2+\alpha_3-
x^*\in \sigma(B_{23}).$ We define $x_i=\lambda_{n-5+i}(B_{23})$,
$i=1,\dots,5.$
 Since
$\rank P_1\le n-9$, then $\lambda_{n-4}(B_{23})\ge \gamma_1$. So
due to assumption, $x_i\notin\{ \alpha_2,\alpha_3,
\alpha_2+\alpha_3\}$ and $\alpha_2+\alpha_3- x_i\in
\sigma(B_{23})$ for every $i=1,\dots,5$. Also $\alpha_2+\alpha_3-
x_1\ge \dots \ge \alpha_2+\alpha_3- x_5$, whence
$\alpha_2+\alpha_3- x_1\ge \lambda_5(B_{23})$.


From inequalities
\begin{align} A-\alpha_1P_1&\ge \diag{0_4,\gamma_1
I_{n-4}}-\alpha_1P_1=\gamma_1 I_n -\alpha_1P_1- \diag{\gamma_1
I_{4},0_{n-4}}\notag \\&\ge (\gamma_1 -\alpha_1)I_n -
\diag{\gamma_1 I_{4},0_{n-4}} = \diag{-\alpha_1 I_{4},(\gamma_1
-\alpha_1)I_{n-4}}.\notag
\end{align}
 we have $\lambda_5(B_{23})\ge \gamma_1-\alpha_1>0$. Now we can
estimate $\alpha_i$: $x_1\ge \gamma_1$ and $
\alpha_2+\alpha_3-x_1\ge \gamma_1-\alpha_1$ $\Longrightarrow$ $
\alpha_1+\alpha_2+\alpha_3 \ge 2\gamma_1$. On the other hand,
 $A\le \gamma_4 I$, so $B_{23}\le \gamma_4 I
-\alpha_1P_1$. By monotonicy property $\lambda_5(B_{23})\le
\lambda_5(\gamma_4 I -\alpha_1P_1)=\gamma_4-\alpha_1$. The
eigenvalue $x_5$ is the biggest eigenvalue of $B_{23}$. By
Proposition \ref{prop1}, $\alpha_2+\alpha_3-x_5$ is the smallest
positive eigenvalue of $B_{23}$, that is $\alpha_2+\alpha_3-x_5
\le \lambda_5(B_{23})\le \gamma_4-\alpha_1$. Taking into account
$x_5\le \gamma_4$, we obtain:
\begin{equation}\label{eq22}
2\gamma_1 \le\alpha_1+\alpha_2+\alpha_3 \le 2\gamma_4.
\end{equation}
Since $\alpha_1> 0$ and it is minimal element of
$\{\alpha_1,\alpha_2,\alpha_3\}$, we get $\alpha_1\le
2\gamma_4/3\approx 2/3$. So $B_{23}$ is invertible and by
Proposition \ref{prop2},
\begin{equation}\label{eq224}
(\alpha_2+\alpha_3-\gamma_4)I\le B_{23}. \end{equation}
 In
addition to this the inequality $\lambda_1(B_{23})\le
\lambda_1(A)=\mu_ 1$ implies $\mu_1\ge
\alpha_2+\alpha_3-\gamma_4$, i.e. applying  left part of
inequalities (\ref{eq22}), $\mu_1\ge
\gamma_1-\alpha_1+(\gamma_1-\gamma_4$). Thus,
\begin{equation}\label{eq221}\alpha_1\ge
(\gamma_1-\mu_1)+(\gamma_1-\gamma_4)\ge \epsilon-\triangle
\end{equation}
 and due to right part of inequalities (\ref{eq22}), we have
\begin{equation}\label{eq222}
\alpha_2\le \gamma_4-\alpha_1/2\le \gamma_4-\epsilon/2 +\triangle
\le \mu_2-10\triangle.\end{equation}
 We are going to localize eigenvalues of
$B_{23}$ more accurately in order to use the same idea as in the
case 2).

Let $$ K_1=\diag{0_{4},0_{18},(\gamma_2-\gamma_1)I_{18},
(\gamma_3-\gamma_1)I_{18},(\gamma_4-\gamma_1)I_{18}}
$$
and
$$
K_2=\diag{\mu_1-\gamma_1,\mu_2-\gamma_1,\mu_3-\gamma_1,\mu_4-\gamma_1,
0_{n-4},}
$$
 Then
 $A=\gamma_1 I+K_1+K_2$. Note that $\|K_1\|\le \triangle$ and
 $\|K_2\|\le \epsilon +\delta$. Also $K_1\ge 0$ and $K_2\le 0$.
 Let $\hat{B}_1=\gamma_1 I+K_1-\alpha_1P_1$ and
 $\hat{B}_2=\gamma_1 I+K_2-\alpha_1P_1$. Counting multiplicity,
 the spectrum of $\hat{B}_1$ has at
 least $k_1$ points   that are less or equal $\gamma_1-\alpha_1+\triangle$. Since $B_{23}\le \hat{B}_1$, we
  have $\lambda_{k_1}(B_{23})\le \gamma_1-\alpha_1+\triangle$.
  Also Theorem \ref{th5} for the sum $\hat{B_2} +\alpha_1 P_1$
   yields $\mu_2=\lambda_{2}(\hat{B}_2 +\alpha_1 P_1)\le
   \lambda_{k_1+2}(\hat{B}_2)$. Since $\hat{B_2}\le B_{23}$, we get
  $\mu_2\le  \lambda_{k_1+2}(B_{23}).$  From inequalities
  (\ref{eq22}) we get $\gamma_1-\alpha_1- \triangle \le
\alpha_2+\alpha_3-\gamma_4$ and so by (\ref{eq224}),
\begin{equation}\label{eq24}
\gamma_1-\alpha_1 -\triangle\le
 \lambda_1(B_{23})\le \lambda_1(B_{23}) \le \dots
 \lambda_{k_1}(B_{23})\le \gamma_1-\alpha_1 +\triangle .
 \end{equation}

 Let us count the number of different eigenvalues of
 $B_{23}$ in the segment $[\mu_2,1-2\triangle]$.
We denote them by $t_1,\dots, t_s$. The eigenvalue $t_i> \alpha_2$
in view of inequality (\ref{eq222}) and $t_i\ne \alpha_2+\alpha_3$
since $\alpha_2+\alpha_3\ge 4/3$. So $t_i=\alpha_3$ or
$t_i>\alpha_3$ and $\alpha_2+\alpha_3 -t_i\in \sigma(B_{23})$. In
the last case we see that
\[\alpha_2+\alpha_3-t_i< \alpha_2\le \gamma_4-\epsilon/2+\triangle<\mu_2
\]
and \[\alpha_2+\alpha_3-t_i\ge \alpha_2+\alpha_3
-(1+2\triangle)\ge 2\gamma_1-\alpha_1-1+2\triangle \ge
\gamma_1-\alpha_1+10\triangle .
\]
As we showed above the only possible eigenvalue of $B_{23}$ from
the interval $(\gamma_4-\alpha_1+\triangle,\mu_2)$ is
$\lambda_{k_1+1}(B_{23})$, that is
$\lambda_{k_1+1}(B_{23})=\alpha_2+\alpha_3-t_i$. Therefore, $s\le
2$ and the set $\sigma(B_{23})\cap
(\gamma_5-\alpha_1+\triangle,1-2\triangle]$ has at most three
points $t_1^*\le t_2^* \le t_3^*$ with the properties $t_1^*\le
\alpha_2<\mu_2\le t_2^* $ and $t_1^*+ t_3^*=\alpha_2+\alpha_3$.
Hence \[ t_3^*-t_1^*\ge \mu_2-\alpha_2=(\mu_2-\gamma_1)+ (\gamma_1
-\alpha_2)\ge
(-\epsilon/6+\triangle)+((\epsilon/2-2\triangle)=\epsilon/3. \]
From definition of $\mu_i$'s we know that
$10\triangle<\mu_4-\mu_3<\mu_3-\mu_2\le \epsilon/18$. So by
Derichlet principle there exist $r\in \{2,3,4\}$ such that
$|\mu_r-t_i^*|>4\triangle $ for every $i=1,2,3$. It means that
\begin{equation}\label{eq25}
\forall i=1,\dots n, |\mu_r-\lambda_i(B_{23})|>4\triangle.
 \end{equation}

Let $h$ be the eigenvector of $A$ with the eigenvalue $\mu_r$. We
define the vector $v$ by the formula $v=P_1 h$. It is a nonzero
vector because $\mu_r\notin \sigma(B_{23}$. Let $P_v$ be the
orthogonal rank one projection defined by $P_v z= (z,v)v/\|v\|^2$.
The operator $\alpha_1P_v+\alpha_2 P_2+\alpha_3P_3$ has the
eigenvalue $\mu_r$ with the eigenvector $h$ by construction.
Denoting $B_3=\alpha_1P_v+B_{23}$, we have $B_3$ is a rank one
perturbation of $B_{23}$. So by interlace theorem
\begin{equation}\label{eq26}
\lambda_i(B_{23})\le \lambda_i(B_{3})\le
\lambda_{i+1}(B_{23}),\quad \forall i=1,\dots n-1
 \end{equation}
and
\begin{equation}\label{eq27}
\sum\limits_{i=1}^n \lambda_i(B_{23})=\sum\limits_{i=1}^n
\lambda_i(B_{3}) -\alpha_1 . \end{equation} Note, that $B_3\le A$,
hence $\lambda_n(B_{23})\le \gamma_4$. Subtracting one part of
(\ref{eq27}) from another and adding $\gamma_4$, we get
\begin{equation}\label{eq277}
\lambda_1(B_{23})+\sum\limits_{i=2}^n (\lambda_i(B_{23})-
\lambda_{i-1}(B_{3}))-\lambda_n(B_{3}) +\alpha_1+\gamma_4=\gamma_4
\end{equation}
which is equivalent to
\begin{equation}\label{eq28}
 \left[\gamma_4-\lambda_n(B_{3})\right]+\sum\limits_{i=2}^n
\left[\lambda_i(B_{23})- \lambda_{i-1}(B_{3})\right]=
\gamma_4-\alpha_1-\lambda_1(B_{23})).
\end{equation}
We note that every summand of (\ref{eq28}) in brackets is
nonnegative. Also one of eigenvalue of $B_3$, say
$\lambda_{i^*}(B_3)$, coincides with $\mu_2$. From (\ref{eq25}) we
conclude that the corresponding expression in the brackets
$\lambda_{i^*}(B_{23})- \lambda_{i^*-1}(B_{3})$ is greater than
$4\triangle$. So the left part of (\ref{eq28}) is greater
$4\triangle$. Taking into account (\ref{eq24}), we see that the
right part of (\ref{eq28}) is less or equal $2\triangle$. So the
equation (\ref{eq28}) is not valid and therefore $A$ is not a
linear combination of $3$ orthoprojections in this case either.
\qed

\begin{cor}\label{cor6} Let $m\in \mathbb{N}$ and $m\ge 76$.
The matrix $\diag{A,\gamma_4 I_{m-76}}$ is not a real linear
combination of $3$ orthoprojections.
    \end{cor}
{\bf Proof}. 
 It is a direct application of the same
arguments as the arguments to $A$ in Proposition \ref{prop6}. \qed

{\bf Concluding remarks.}

1. The scheme of the proof of Theorem \ref{th1} can be directly
applied to decompositions of finite matrices in unitary space,
since every Hermitian matrix with zero trace is a self-commutator
\cite{Thompson58}. For example, for a $2n\times 2n $ matrix $A$,
we put $\lambda=\tr A/n$ and then take all steps according to the
proof. If $A$ is a $2n+1\times 2n+1 $ matrix, it is enough to
consider the case $A=\diag{\mu,A_1}$ with $|\mu| =\|A\|$. Here the
orthoprojections $P_i$ will be of the form
$\hat{P}_i=\diag{1,{P_i}}$, $i=1,3$ and
$\hat{P}_i=\diag{0,{P_i}}$, $i=2,4$ where $P_i$ are
orthoprojections from the decomposition of $A_1$ into the linear
combinations of $4$ orthoprojections with the restriction
$a-c=\mu$ on the coefficients in the proof of Theorem \ref{th1} .

2. In view of Proposition \ref{prop6}, it is interesting to know
what is the maximal number $n$ for which every Hermitian $k\times
k$ matrix is a real linear combination of three orthoprojections
providing $k\le n$. We suppose it is not greater than $25$ since
many cases of the proof of Proposition \ref{prop6} can be applied
directly for smaller value of $n$.

   INSTITUTE OF MATHEMATICS,
    Tereshchenkivs'ka 3, Kyiv, 01601,  Ukraine\\
  \emph{E-mail address}: slavik@imath.kiev.ua

\end{document}